\documentclass[twocolumn]{svjour3}
\smartqed
\usepackage{epsfig,amsmath,amsfonts,amssymb,color,hyperref}%
\usepackage{subcaption}
\usepackage{tikz}
\usepackage{pgfplots}
\newcommand{\curl}{\operatorname{curl}}

\newcommand{\isdef}{\mathrel{\mathrel{\mathop:}=}}

\renewcommand{\d}{\operatorname{d}\!}
\newcommand{\bs}{\boldsymbol}

\begin{document}
\title{Solving acoustic scattering problems by the isogeometric boundary element method
\thanks{The first author acknowledges the support by the Deutsche 
Forschungsgemeinschaft (DFG, German Research Foundation) 
under Germany’s Excellence Strategy -- GZ 2047/1, Project-ID 390685813.}}
\author{J\"urgen D\"olz \and Helmut Harbrecht \and Michael Multerer}
\institute{J\"urgen D\"olz \at
Institut f\"ur Numerische Simulation,
Universit\"at Bonn, Friedrich-Hirzebruch-Allee 7, 53115 Bonn, Germany\\
\email{doelz@ins.uni-bonn.de}
\and Helmut Harbrecht \at
Departement Mathematik und Informatik,
Universit\"at Basel, Spiegelgasse 1, 4051 Basel, Schweiz\\
\email{helmut.harbrecht@unibas.ch}
\and Michael Multerer \at
Istituto Eulero, Universit\`a della Svizzera italiana,
Via la Santa 1, 6962 Lugano, Svizzera\\
\email{michael.multerer@usi.ch}}
\date{Received: date / Accepted: date}
\maketitle
\begin{abstract}
We solve acoustic scattering problems by means of the 
isogeometric boundary integral equation method. In order to avoid
spurious modes, we apply the combined field integral equations
for either sound-hard scatterers or sound-soft scatterers. These 
integral equations are discretized by Galerkin's method, which 
especially enables the mathematically correct regularization of 
the hypersingular integral operator. In order to circumvent densely 
populated system matrices, we employ the isogeometric fast multipole
method. The result is an algorithm that scales essentially linear
in the number of boundary elements. Numerical experiments are
performed which show the feasibility and the performance of the
approach.
\end{abstract}

\section{Introduction}
Acoustic wave scattering appears in many places in engineering
practice. This includes, for instance, the modelling of sonar and 
other methods of acoustic location, as well as outdoor noise 
propagation and control, especially stemming from automobiles, 
railways or aircrafts. Since an analytical solution of scattering 
problems is in general impossible, numerical approaches are 
called for the approximate solution.

Most acoustic scattering problems may be formulated in the 
frequency domain by employing the Helm\-holtz equation. Assume 
that an acoustic wave encounters an impenetrable, bounded obstacle 
$\Omega\subset\mathbb{R}^3$, having a Lipschitz smooth boundary $\Gamma
\isdef\partial\Omega$, and, as a consequence, gets scattered. Given
the \emph{incident plane wave} $u_{\text{inc}}({\bs x}) = e^{i\kappa\langle{\bs d},
{\bs x}\rangle}$ with known wavenumber $\kappa$ and direction ${\bs d}$,
where $\|{\bs d}\|_2=1$, the goal is to compute the \emph{scattered wave}
$u_{\mathrm{s}}$. The physical model behind this is as follows. The 
\emph{total wave} $u = u_{\text{inc}}+u_{\text{s}}$ satisfies the 
exterior Helmholtz equation
\begin{equation}\label{eq:pde1}
  \Delta u + \kappa^2 u = 0\ \text{in}\ \mathbb{R}^3\setminus\overline{\Omega}.
\end{equation}
The boundary condition at the scatterer's surface
depends on its physical properties. If the scatterer 
constitutes a \emph{sound-soft} obstacle, then 
the acoustic pressure vanishes at \(\Gamma\)
and we
have the homogeneous Dirichlet condition 
\begin{equation}\label{eq:pde2a}
  u = 0\ \text{on}\ \Gamma.
\end{equation}
Whereas, if the scatterer constitutes a \emph{sound-hard} 
obstacle, then the normal velocity vanishes at \(\Gamma\)
and we have the homogeneous Neumann condition 
\begin{equation}\label{eq:pde2b}
  \frac{\partial u}{\partial \bs n} = 0\ \text{on}\ \Gamma.
\end{equation}
The behaviour towards infinity is imposed by the
Sommerfeld radiation condition
\begin{equation}	\label{eq:pde3}
  \lim_{r\to\infty} r \left\{ \frac{\partial u_s}{\partial r} 
  	- i\kappa u_s\right\} = 0,\ \text{where}\ r\isdef\|\bs x\|_2.
\end{equation}
It implies the asymptotic expansion
\[
  u_s(\bs x) = \frac{e^{i\kappa\|\bs x\|_2}}{\|\bs x\|_2} 
  	\left\{u_\infty\Big(\frac{\bs x}{\|\bs x\|_2}\Big)
	+ \mathcal{O}\Big(\frac{1}{\|\bs x\|_2}\Big)\right\}
\]
as $\|\bf x\|_2\to\infty$. Herein, the function $u_\infty\colon\mathbb{S}^1
:=\{\hat{\bs x}\in\mathbb{R}^d:\|\hat{\bs x}\|_2=1\}\to\mathbb{C}$
is called the \emph{far-field pattern}, which is always analytic 
in accordance with \cite[Chapter 6]{CK2}. In applications, the 
far-field pattern is the most important quantity of interest derived
from scattering problems.

To avoid the discretization of the unbounded exterior 
domain $\mathbb{R}^3\setminus\overline{\Omega}$, one can 
exploit the integral equation formalism to compute the numerical 
solution of acoustic scattering problems. Then, one arrives at 
a boundary integral equation only defined on the boundary 
$\Gamma$. We will employ here the methodology of
\emph{isogeometric analysis} (IGA) to discretize this boundary 
integral equation. IGA has been introduced in \cite{Hughes_2005aa} 
in order to incorporate simulation techniques into the design 
workflow of industrial development. The goal is thus to unify 
the CAD representation of the scatterer with the boundary element 
discretization of the integral equation in terms of \emph{non-uniform 
rational B-splines} (NURBS). We refer the reader to 
\cite{DHK+18,MZBF,SBTR12} and the references therein 
for details of the isogeometric boundary element method.

While a reformulation of the scattering problem by means 
of a boundary integral equation replaces the problem posed in 
the unbounded domain by a problem posed on the scatterer's
closed boundary, the linear operator under consideration becomes 
a nonlocal boundary integral operator. This results in densely 
populated matrices with the consequence that desirable 
realistic simulations would still beyond current computing 
capacities. Therefore, we combine the isogeometric boundary
element method with the fast multipole method proposed in
\cite{GR}. Our particular implementation relies on interpolation 
as proposed in \cite{HP13} and yields a black-box version 
of the fast multipole method which is applicable to any
asymptotically smooth integral kernel.

The isogeometric boundary element method we use 
has been developed in \cite{DHK+18,DHP16} and was
made accessible to the public by the software \texttt{C++} library
\verb+Bembel+ \cite{bembel,bembelpaper}. \verb+Bembel+ has 
for example been applied successfully to engineering problems arising 
from electromagnetics \cite{DKSW20,KSUW21} or from acoustics 
\cite{DHJM23}. It has also been used in other applications,
for example, to optimize periodic structures \cite{HMR22}, 
in uncertainty quantification \cite{Doe20,DHJM23}, the
coupling of FEM and BEM \cite{EEK21}, or the partial
element equivalent circuit (PEEC) method \cite{TNSR23}. 
In this article, we shall present results for the frequency stable
solution of acoustic obstacle scattering problems by using 
combined field integral equations. Although we restrict ourselves
to the sound-soft and sound-hard cases, the presented concepts 
are also suitable to treat penetrable obstacles, i.e.~objects 
described by a different diffractive index to the free space.
 
The rest of this article is structured as follows.
In Section~\ref{sec:implementation}, we introduce the
frequency stable boundary integral equations which are employed
to solve either sound-hard or sound-soft scattering problems. 
Section~\ref{sec:IGA} recapitulates the basic concepts from
isogemeometric analysis and introduces the discretization spaces
that will be used later on. In Section~\ref{sec:discrete},
we discuss the discretization of the required boundary integral
operators. In particular, we address the regularization of 
the hypersingular operator. Moreover, we comment on the
isogeometric fast multipole method for the fast assembly of
the operators and the potential evaluation. The numerical
experiments are presented in Section~\ref{sec:experiments},
while concluding remarks are stated in Section~\ref{sec:conclusion}.
 
\section{Boundary integral equation method}\label{sec:implementation}
In order to solve the boundary value problem
\eqref{eq:pde1}--\eqref{eq:pde3}, we shall employ 
a suitable reformulation by boundary integral equations.
To this end, we introduce the acoustic single layer operator
\begin{align*}
  &\mathcal{V}\colon H^{-1/2}(\Gamma)\to H^{1/2}(\Gamma),\\
  &\qquad(\mathcal{V}\rho)({\bs x}):=\int_{\Gamma}
  	G({\bs x},{\bs y})\rho({\bs y})\d\sigma_{\bs y},
\end{align*}
the acoustic double layer operator
\begin{align*}
  &\mathcal{K}\colon L^2(\Gamma)\to L^2(\Gamma),\\
  &\qquad(\mathcal{K}\rho)({\bs x}):=\int_{\Gamma}
  	\frac{\partial G({\bs x},{\bs y})}
	{\partial{\bs n}_{\bs y}}\rho({\bs y})\d\sigma_{\bs y},
\end{align*}
its adjoint
\begin{align*}
  &\mathcal{K}^\star\colon L^2(\Gamma)\to L^2(\Gamma),\\
  &\qquad(\mathcal{K}^\star\rho)({\bs x}):=\int_{\Gamma}
  	\frac{\partial G({\bs x},{\bs y})}
	{\partial{\bs n}_{\bs x}}\rho({\bs y})\d\sigma_{\bf y},
\end{align*}
as well as the acoustic hypersingular operator
\begin{equation}\label{eq:hypersingular}
\begin{aligned}
  &\mathcal{W}\colon H^{1/2}(\Gamma)\to H^{-1/2}(\Gamma),\\
  &\qquad(\mathcal{W}\rho)({\bs x}):=-\frac{1}{\partial{\bs n}_{\bs x}}\int_{\Gamma}
  	\frac{\partial G({\bs x},{\bs y})}{\partial{\bs n}_{\bs y}}
		\rho({\bs y})\d\sigma_{\bs y}.
\end{aligned}
\end{equation}
Here, ${\bs n}_{\bs x}$ and ${\bs n}_{\bs y}$ denote the 
outward pointing normal vectors at the surface points 
${\bs x}, {\bs y}\in\Gamma$, respectively, while $G(\cdot,\cdot)$ 
denotes the fundamental solution for the Helmholtz equation. In three 
spatial dimensions, the latter is given by
\[
  G({\bs x}, {\bs y}) = \frac{e^{i\kappa\|{\bs x}-{\bs y}\|_2}}{4\pi\|{\bs x}-{\bs y}\|_2}.
\]

Although the Helmholtz problem \eqref{eq:pde1}--\eqref{eq:pde3}
is uniquely solvable, a respective boundary integral formulation might
not if $\kappa^2$ is an eigenvalue for the Laplacian inside the scatterer
$\Omega$. In order to avoid such \emph{spurious modes}, we employ 
combined field integral equations in the following. Then, for some 
real $\eta\neq 0$, the solution of the boundary integral equation
\begin{equation}\label{eq:D2N}
  \bigg(\frac{1}{2}+\mathcal{K}^\star-i\eta\mathcal{V}\bigg)
  	\frac{\partial u}{\partial{\bs n}} = \frac{\partial u_{\text{inc}}}{\partial{\bs n}}
		-i\eta u_{\text{inc}}
\end{equation}
gives rise to the scattered wave in accordance with
\begin{equation}\label{eq:pot1}
  u_{\text{s}}({\bs x}) = \int_\Gamma G({\bs x},{\bs y})
  \frac{\partial u({\bs y})}{\partial{\bs n}_{\bs y}}\d\sigma_{\bs y}
\end{equation}
in case of sound-soft scattering problems.
In case of sound-hard obstacles, we will solve the integral equation
\begin{equation}\label{eq:N2D}
  \bigg(\frac{1}{2} - \mathcal{K} + i\eta \mathcal{W}\bigg)u 
  = u_{\text{inc}}-i\eta\frac{\partial u_{\text{inc}}}{\partial\bs n}.
\end{equation}
Having solved \eqref{eq:N2D}, the scattered wave 
is computed by
\begin{equation}\label{eq:pot2}
  u_{\text{s}}({\bs x}) = \int_\Gamma 
  \frac{\partial G({\bs x},{\bs y})}{\partial{\bs n}_{\bs y}}u({\bs y})\d\sigma_{\bs y}.
\end{equation}
Notice that the boundary integral equations \eqref{eq:D2N} and
\eqref{eq:N2D} are always uniquely solvable, independent of 
the wavenumber $\kappa$, compare \cite{BM,CK2,Kuss1969}.

\section{Isogeometric analysis}\label{sec:IGA}
\subsection{B-splines}
We shall give a brief introduction to the basic concepts 
of isogeometric analysis, starting with the definition of the 
B-spline basis, followed by the description of the scatterer 
by using NURBS. To this end, let $\mathbb{K}$ be either 
$\mathbb{R}$ or $\mathbb{C}$. The original definitions (or 
equivalent notions) and proofs, as well as basic algorithms, 
can be found in most of the standard spline and isogeometric 
literature \cite{Cottrell_2009aa,Hughes_2005aa,%
Piegl_1997aa,Schumaker_1981aa,Lee_1996aa}.

\begin{figure*}
\begin{subfigure}{.5\textwidth}
\begin{tikzpicture}
\begin{axis}[
xmin = -.1,
xmax = 1.1,
ymin = -.1,
ymax = 1.1,
width=1\columnwidth,
height=.66\columnwidth,
grid=major,
legend style={
at={(.5,1)},
anchor=south}
]
\addplot[red,ultra thick,mark=none,domain = 0:1/3] {1};
\addplot[teal,ultra thick,mark=none,domain = 1/3:2/3] {1};
\addplot[blue,ultra thick,mark=none,domain = 2/3:1] {1};
\end{axis}
\end{tikzpicture}
\caption{$p=0$, $\Xi=[0,1/3,2/3,1]$.}
\end{subfigure}\!\!
\begin{subfigure}{.5\textwidth}
\begin{tikzpicture}
\begin{axis}[
xmin = -.1,
xmax = 1.1,
ymin = -.1,
ymax = 1.1,
width=1\columnwidth,
height=.66\columnwidth,
grid=major,
legend style={
at={(.5,1)},
anchor=south}
]
\addplot[red,ultra thick,mark=none,domain = 0:1/3] {3*-(x-1/3)};
\addplot[teal,ultra thick,mark=none,domain = 0:1/3] {3*(x)};
\addplot[blue,ultra thick,mark=none,domain = 1/3:2/3] {3*(x-1/3)};
\addplot[orange,ultra thick,mark=none,domain = 2/3:1]{3*(x-2/3)};
\addplot[teal,ultra thick,mark=none,domain = 1/3:2/3] {(1-3*x)+1};
\addplot[blue,ultra thick,mark=none,domain = 2/3:1] {3*(1-x-1/3)+1};
\end{axis}
\end{tikzpicture}
\caption{$p=1$, $\Xi=[0,0,1/3,2/3,1,1]$.}
\end{subfigure}\\[.5cm]
\begin{subfigure}{\textwidth}
\begin{tikzpicture}
\begin{axis}[
xmin = -.1,
xmax = 1.1,
ymin = -.1,
ymax = 1.1,
width=1\columnwidth,
height=.3397\columnwidth,
grid=major,
legend style={
at={(1,.5)},
anchor=west}
]
\addplot[red,ultra thick,mark=none,domain = 0:1/3]{(3*(x)-1)^2} ;
\addplot[teal,ultra thick,mark=none,domain = 0:1/3]{2*(3*(x))*(1-3*x)+.5*(3*x)^2};
\addplot[blue,ultra thick,mark=none,domain = 0:1/3]{.5*(x*3)^2};
\addplot[orange,ultra thick,mark=none,domain = 2/3:1]{2*(3*((x-2/3)))*(1-3*(x-2/3))+.5*(1-3*(x-2/3))^2};
\addplot[brown,ultra thick,mark=none,domain = 2/3:1]{(3*(x-2/3))^2};
\addplot[teal,ultra thick,mark=none,domain = 1/3:2/3]{.5*(1-3*(x-1/3))^2 };
\addplot[blue,ultra thick,mark=none,domain = 1/3:2/3]{-(((x*3)-1)-1)*((x*3)-1)+.5};
\addplot[blue,ultra thick,mark=none,domain = 2/3:1]{.5*(3-(x*3))^2};
\addplot[orange,ultra thick,mark=none,domain = 1/3:2/3]{.5*(3*(x-1/3))^2 };
\end{axis}
\end{tikzpicture}
\caption{$p=2$, $\Xi=[0,0,0,1/3,2/3,1,1,1]$.}
\end{subfigure}
\caption{B-spline bases for $p=0,1,2$ and open knot 
vectors with interior knots $1/3$ and $2/3$.}\label{fig::splines}
\end{figure*}
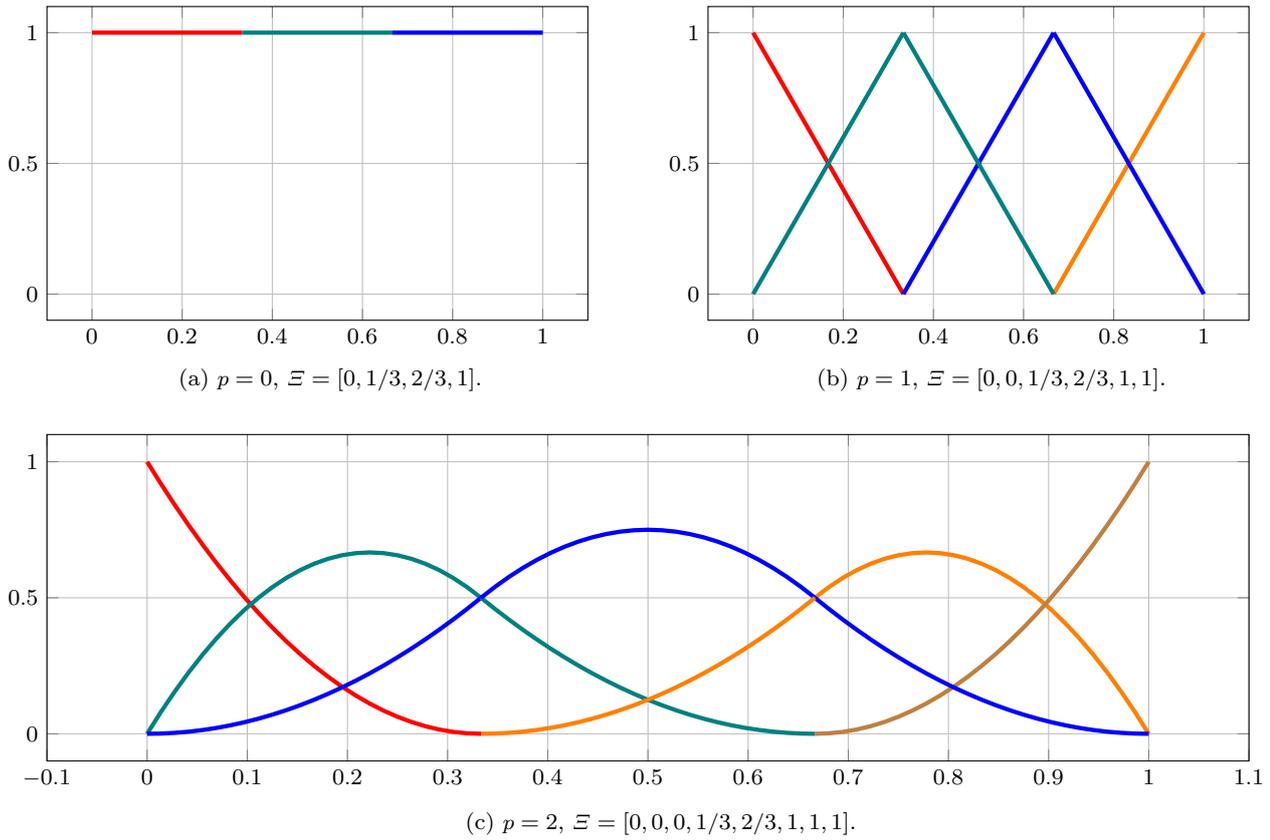

\begin{definition}
Let $0\leq p\leq k$. We define a \emph{$p$-open knot 
vector} as a set
\begin{align*}
  &\Xi = \big[\underbrace{\xi_0 = \cdots =\xi_{p}}_{=0}\leq\cdots\\ 
  &\hspace*{2cm}\cdots\leq \underbrace{\xi_{k}=\cdots =\xi_{k+p}}_{=1}\big]
  \in[0,1]^{k+p+1},
\end{align*}
where $k$ denotes the number of control points.
The associated basis functions are given by
$\{b_j^p\}_{j=0}^{k-1}$ for $p=0$ as
\[
  b_j^0(x) =\begin{cases}
	1, & \text{if }\xi_j\leq x<\xi_{j+1}, \\
	0, & \text{otherwise},\end{cases}
\]
and for $p>0$ via the recursive relationship
\[
  b_j^p(x) = \frac{x-\xi_j}{\xi_{j+p}-\xi_j}b_j^{p-1}(x) 
  +\frac{\xi_{j+p+1}-x}{\xi_{j+p+1}-\xi_{j+1}}b_{j+1}^{p-1}(x),
\]
cf.~Figure \ref{fig::splines}. A \emph{B-spline} is then 
defined as a function
\[
f(x) = \sum_{0\leq j< k}p_jb_j^p(x),
\]
where $\{p_j\}_{j=0}^{k-1}\subset\mathbb{K}$ 
denotes the set of \emph{control points}. If one 
sets $\{p_j\}_{j=0}^{k-1}\subset\mathbb{K}^d$, 
then $f$ will be called a \emph{B-spline curve}.
\end{definition}

Having the spline functions at hand, we can introduce
the spline spaces which serve as fundament for the
definition of the ansatz and test spaces of the boundary
element method.

\begin{definition}
Let $\Xi$ be a $p$-open knot vector containing $k+p+1$ 
elements. We define the \emph{spline space} $S_{p}(\Xi)$ 
as the space spanned by $\{b_j^p\}_{j=0}^{k-1}$.
\end{definition}

Finally, we should consider the relation between the 
spline spaces and the underlying mesh relative to a 
certain mesh size.

\begin{definition}
For a knot vector $\Xi,$ we define the \emph{mesh size} 
$h$ to be the maximal distance 
\begin{equation}\label{eq:h}
h\isdef \max_{j=0}^{k+p-1}h_j,\quad\text{where}\quad
	h_j\isdef\xi_{j+1}-\xi_{j},
\end{equation}
between neighbouring knots.
We call a knot vector \emph{quasi uniform}, when there 
exists a constant $\theta\geq 1$ such that for all $j$ the 
ratio $h_j\cdot h_{j+1}^{-1}$ satisfies $\theta^{-1}\leq 
h_j\cdot h_{j+1}^{-1} \leq \theta.$
\end{definition}

B-splines on higher dimensional domains are 
constructed through simple tensor product relationships 
for $p_{j_1,\dots j_\ell}\in\mathbb{K}^d$ via
\begin{equation}\label{def::tpspline}
\begin{aligned}
  &f(x_1,\dots ,x_\ell)\\
  &\quad =\sum_{j_1=0}^{k_1-1}\dots
  	\sum_{j_\ell=0}^{k_\ell-1} p_{j_1,\dots,j_\ell} \cdot 
	b_{j_1}^{p_1}(x_1)\cdots b_{j_\ell}^{p_\ell}(x_\ell),
\end{aligned}
\end{equation}
which allows \emph{tensor product B-spline spaces} to 
be defined as 
\[
S_{p_1,\dots,p_\ell}(\Xi_1,\dots,\Xi_\ell).
\]

Throughout this article, we will reserve the letter $h$ for 
the mesh size \eqref{eq:h}. All knot vectors will be assumed to 
be quasi uniform, such that the usual spline theory is applicable 
\cite{BuffaActa,Piegl_1997aa,Schumaker_1981aa}.

\subsection{Isogeometric representation of the scatterer}
We assume that the boundary $\Gamma$ of the scatterer 
is closed and Lipschitz continuous. For the remainder
of this article, we assume that it is given patchwise as
$\Gamma=\bigcup_{j=1}^n\Gamma_j$, i.e.~that it is
induced by (smooth, nonsingular, bijective) diffeomorphisms
\begin{equation}\label{def::geom}
  {\bf F}_j\colon \widehat \Omega = [0,1]^2 \to \Gamma_j \subset \mathbb{R}^3.
\end{equation}
In the spirit of isogemetric analysis, these mappings 
are given by NURBS mappings, i.e.~by
\[
  {\bf F}_j(x,y)\isdef \sum_{j_1=0}^{k_1-1}\sum_{j_2=0}^{k_2-1}
  \frac{c_{j_1,j_2} b_{j_1}^{p_1}(x) b_{j_2}^{p_2}(y) w_{j_1,j_2}}
  {\sum_{i_1=0}^{k_1-1}\sum_{i_2=0}^{k_2-1} b_{i_1}^{p_1}(x) b_{i_2}^{p_2}(y) w_{i_1,i_2}}
\]
with control points $c_{j_1,j_2}\in \mathbb{R}^3$ and weights 
$w_{i_1,i_2}>0$. We will moreover require that, for any interface 
$D = \Gamma_j\cap \Gamma_i \neq \emptyset$, the NURBS 
mappings coincide, i.e.~that, up to rotation of the reference 
domain, one finds ${\bf F}_j(\cdot,1) \equiv {\bf F}_i(\cdot,0)$.

\subsection{Ansatz and test spaces}
The mappings of \eqref{def::geom} give rise to the transformations
\[
  \iota_j(f) \isdef f\circ {\bf F}_j,
\]
which can be utilized to define discrete spaces patchwise, 
by mapping the space of tensor product B-splines as in 
\eqref{def::tpspline} with
\[
  \Xi_{p,m} \isdef \big[ \underbrace{0,\dots,0}_{p+1\text{ times}},
  	{1}/{2^m},\dots,  {(2^m-1)}/{2^m},
	\underbrace{1,\dots,1}_{p+1\text{ times}}\big]
\]
to the geometry. Here, the variable $m$ denotes the level 
of uniform refinement. For the purposes of discretizing $\mathcal{V}$, $\mathcal{K}$, and $\mathcal{K}^\star$, the global function 
space on $\Gamma$ defined by
\begin{align*}
  &\mathbb S_{p,m}^2(\Gamma) \isdef \Big\{
  	f\in H^{-1/2}(\Gamma)\colon f_{|\Gamma_j} \equiv \iota_j^{-1}(g)\\
	&\hspace*{2cm}\text{ for some }g\in S_{p,p}(\Xi_{p,m},\Xi_{p,m})\Big\},
\end{align*}
as commonly done in the isogeometric literature, see 
e.g.~\cite{Buffa_2011aa,Buffa_2020aa}, is sufficient. 
Note that the spline space $\mathbb S_{p,m}^2(\Gamma)$ 
is of dimension $n\cdot(2^m + p )^2$, where $n$ denotes 
the number of patches involved in the description of the geometry. For the purposes of discretizing $\mathcal{W}$, we also require the space 
\begin{equation}\label{def::space2}
\begin{aligned}
&\mathbb S_{p,m}^0(\Gamma) \isdef \Big\{
f\in H^{1/2}(\Gamma)\colon f_{|\Gamma_j} \equiv \iota_j^{-1}(g)\\
&\hspace*{2cm}\text{ for some }g\in S_{p,p}(\Xi_{p,m},\Xi_{p,m})\Big\},
\end{aligned}
\end{equation}
see, e.g., also \cite{Buffa_2011aa,Buffa_2020aa}. Note that that 
$\mathbb S_{p,m}^0(\Gamma)\subset \mathbb S_{p,m}^2(\Gamma)$ 
consists of globally continuous B-splines whereas $\mathbb S_{p,m}^2(\Gamma)$ 
is discontinuous across patch boundaries.

\section{Discretization}\label{sec:discrete}
\subsection{Galerkin method}
With the boundary integral equations and a collection of 
spline spaces available, we are now in the position to discretize 
\eqref{eq:D2N} and \eqref{eq:N2D}. We consider a Galerkin 
discretization in the $L^2(\Gamma)$-duality product with 
the spline spaces $\mathbb S_{p,m}^2(\Gamma)$ and 
$\mathbb S_{p,m}^0(\Gamma)$ as ansatz and test spaces. Thus, 
the discrete variational formulation for \eqref{eq:D2N} reads
\begin{align*}
&\text{Find}~t_h\in \mathbb S_{p,m}^2(\Gamma)~\text{such that}\\
&\hspace*{1cm}
\frac{1}{2}\langle t_h,v_h\rangle_{\Gamma}+\langle\mathcal{K}^\star t_h,v_h\rangle_{\Gamma}-i\eta\langle\mathcal{V}t_h,v_h\rangle_{\Gamma}\\
&\hspace*{2cm}=\Big\langle \frac{\partial u_{\text{inc}}}{\partial{\bs n}}
		-i\eta u_{\text{inc}},v_h\Big\rangle_{\Gamma}\\
&\text{for all}~v_h\in \mathbb S_{p,m}^2(\Gamma),
\end{align*}
with the Galerkin approximation $t_h\approx\partial u/\partial{\bf n}$. 
Choosing a basis $\mathbb S_{p,m}^2(\Gamma)=\operatorname{span}
\{\psi_{2,1},\ldots,\psi_{2,N}\}$ leads to the system of linear equations
\begin{equation}\label{eq:CFIED2N}
\bigg(\frac{1}{2}{\bf M}_2+{\bf K}_2^\star-i\eta{\bf V}_2\bigg){\bf t}={\bf u}_2
\end{equation}
with
\begin{align*}
{\bf M}_2&=\big[\langle\psi_{2,j},\psi_{2,i}\rangle_{\Gamma}\big]_{i,j=1}^N,\\
{\bf K}_2^{\star}&=\big[\langle\mathcal{K}^\star\psi_{2,j},\psi_{2,i}\rangle_{\Gamma}\big]_{i,j=1}^N,\\
{\bf V}_2&=\big[\langle\mathcal{V}\psi_{2,j},\psi_{2,i}\rangle_{\Gamma}\big]_{i,j=1}^N,\\
{\bf u}_2&=\Big[\Big\langle \frac{\partial u_{\text{inc}}}{\partial{\bs n}}
		-i\eta u_{\text{inc}},\psi_{2,i}\Big\rangle_{\Gamma}\Big]_{i=0}^N,
\end{align*}
and ${\bf t}$ being the coefficient vector of $t_h$.

The discrete variational formulation for \eqref{eq:N2D} reads
\begin{align*}
&\text{Find}~g_h\in \mathbb S_{p,m}^0(\Gamma)~\text{such that}\\
&\hspace*{1cm}
\frac{1}{2}\langle g_h,v_h\rangle_{\Gamma}+\langle\mathcal{K} g_h,v_h\rangle_{\Gamma}-i\eta\langle\mathcal{W}g_h,v_h\rangle_{\Gamma}\\
&\hspace*{2cm}=\Big\langle u_{\text{inc}}-i\eta\frac{\partial u_{\text{inc}}}{\partial\bs n},v_h\Big\rangle_{\Gamma}\\
&\text{for all}~v_h\in \mathbb S_{p,m}^0(\Gamma),
\end{align*}
with the Galerkin approximation $g_h\approx u|_\Gamma$. Choosing 
a basis $\mathbb S_{p,m}^0(\Gamma)=\operatorname{span}\{\psi_{0,1},
\ldots,\psi_{0,M}\}$ leads to the linear system of equations
\begin{equation}\label{eq:CFIEN2D}
\bigg(\frac{1}{2}{\bf M}_0-{\bf K}_0-i\eta{\bf W}_0\bigg){\bf g}={\bf v}_0
\end{equation}
with
\begin{align*}
{\bf M}_0&=\big[\langle\psi_{0,j},\psi_{0,i}\rangle_{\Gamma}\big]_{i,j=1}^M,\\
{\bf K}_0&=\big[\langle\mathcal{K}\psi_{0,j},\psi_{0,i}\rangle_{\Gamma}\big]_{i,j=1}^M,\\
{\bf W}_0&=\big[\langle\mathcal{W}\psi_{0,j},\psi_{0,i}\rangle_{\Gamma}\big]_{i,j=1}^M,\\
{\bf v}_0&=\bigg[\bigg\langle u_{\text{inc}}-i\eta
\frac{\partial u_{\text{inc}}}{\partial\bs n},\psi_{0,i}\bigg\rangle_{\Gamma}\bigg]_{i=0}^M,
\end{align*}
and ${\bf g}$ being the coefficient vector of $g_h$.

It is well known that the matrices ${\bf V}_2$, ${\bf K}_0$, 
${\bf K}_2^\star$, and ${\bf W}_0$ are dense, which makes 
the assembly and storage of these matrices as well as the 
solution of the corresponding linear systems of equations 
computationally prohibitively expensive for higher resolution 
of the ansatz spaces, i.e., large $M$ or $N$. This is why
we shall apply the multipole method presented in Subsection~\ref{sec:multipole}.

\subsection{Reformulation on the reference domain}\label{sec:bilinearonreference}
Due to the isogeometric representations of the geometry, the 
bilinear forms for the computation of the matrix entries can entirely 
pulled back to the reference domain \cite{HP13}. To this end, 
let $\mathcal{A}$ with
\begin{equation}\label{eq:integraloperator}
(\mathcal{A}\mu)({\bs x})= \int_\Gamma k({\bs x},{\bs y})\mu({\bs y})\d\sigma_{\bs y},\qquad {\bs x}\in\Gamma
\end{equation}
be one of the operators $\mathcal{V}$, $\mathcal{K}$, or 
$\mathcal{K}^\star$ and $\mu,\nu\colon\Gamma\to\mathbb{C}$ 
be functions of sufficient regularity. Defining the \emph{surface 
measure} of a mapping ${\bf F}_j$ for $\hat{\bs x} = (x,y)\in [0,1]^2$ as
\[
  a_j (\hat{\bs x})\isdef\big\|\partial_{x}{\bf F}_j(\hat{\bs x})
  	\times \partial_{y}{\bf F}_j(\hat{\bs x})\big\|_2,
\]
the bilinear forms for the matrix entries can be recast as
\begin{align*}
  &\langle\mathcal{A}\mu,\nu\rangle_\Gamma
  = \sum_{j=1}^n \langle\mathcal{A}\mu,\nu\rangle_{\Gamma_j}\\
  &\quad= \sum_{i,j=1}^n \int_{\Gamma_i}\int_{\Gamma_j}  
  	k({\bf x},{\bf y}) \mu({\bs x})\nu({\bs y})\d\sigma_{\bs y}\d\sigma_{\bs x}\\
  &\quad= \sum_{i,j=1}^n\int_{[0,1]^2}\int_{[0,1]^2}  
  	k\big({\bf F}_j(\hat{\bs x}),{\bf F}_i(\hat{\bs y})\big)\\
  &\hspace*{2cm}\times\mu\big({\bf F}_j(\hat{\bs x})\big)
                 \nu\big({\bf F}_i(\hat{\bs y})\big)
		a_{j}(\hat{\bs x})a_{i}(\hat{\bs y})
		\d\hat{\bs y}\d\hat{\bs x}\\
  &\quad= \sum_{i,j=1}^n\int_{[0,1]^2}\int_{[0,1]^2} 
  	k_{j,i}(\hat{\bs x},\hat{\bs y})\mu_j(\hat{\bs x})\nu_i(\hat{\bs y})\d\hat{\bs y}\d\hat{\bs x},
\end{align*}
with the pull-back of the kernel function and the ansatz 
and test functions
\begin{equation}\label{eq:pbkernel}
\begin{aligned}
k_{j,i}(\hat{\bs x},\hat{\bs y})&=a_{j}(\hat{\bs x})a_{i}(\hat{\bs y})k\big({\bf F}_j(\hat{\bs x}),{\bf F}_i(\hat{\bs y})\big),\\
\mu_j(\hat{\bs x})&=\iota_j(\mu)(\hat{\bs x}),\\
\nu_i(\hat{\bs y})&=\iota_i(\nu)(\hat{\bs y}).
\end{aligned}
\end{equation}
Applying a similar reasoning to the right-hand 
side yields
\[
  \langle g,\nu\rangle_\Gamma 
  	= \sum_{i=1}^n \int_{[0,1]^2} 
	g\big({\bf F}_i(\hat{\bs x})\big)\nu_i(\hat{\bs x})a_{i}(\hat{\bs x})\d\hat{\bs x}.
\]
Due to the additional derivative, the hypersingular operator 
$\mathcal{W}$ requires a special treatment which we will 
elaborate next.

\subsection{Regularization of the Helmholtz hypersingular operator}
\label{sec:regularizedhyper}
The hypersingular operator $\mathcal{W}$ from
\eqref{eq:hypersingular} does not have a well defined integral operator
representation as in \eqref{eq:integraloperator}. Instead, it is common
knowledge that the operator can be replaced by a regularized one
in case of a Galerkin discretization. Namely, for the computation 
of the matrix entries, the representation
\[
\begin{aligned}
&\langle\mathcal{W}\psi_{0,j},\psi_{0,i}\rangle_{\Gamma} 
=\langle\mathcal{V}\curl_{\Gamma}\psi_{0,j},\curl_{\Gamma}\psi_{0,i}\rangle_{\Gamma}\\
&\ -
\kappa^2\int_{\Gamma}\int_{\Gamma}G({\bs x},{\bs y})
\langle{\bs n}_{\bs x},{\bs n}_{\bs y}\rangle_{\mathbb{R}^3}
\psi_{0,j}({\bs x})\psi_{0,i}({\bs y})\d\sigma_{\bs y}\d\sigma_{\bs x},
\end{aligned}
\]
$i,j=1,\ldots,M$, can be used, cf.\ e.g.~\cite{Maue1949}. Therein,
$\curl_{\Gamma}\psi_{0,i}$ denotes the surface curl which maps 
a scalar valued function on the surface into a vector field in the 
tangential space of $\Gamma$. On any given patch $\Gamma_j$,
the isogeometric representations of the boundary of the scatterer 
allow for its explicit representation
\begin{equation}\label{eq:surfacecurl}
\begin{aligned}
&\curl _{\Gamma}\psi_{0,i}(\mathbf{x})=\frac{1}{a_i(\hat{\mathbf{x}})}
\Big(\partial_{\hat{x}_1}\iota_j(\psi_{0,i})(\hat{\mathbf{x}})
\partial_{\hat{x}_2}{\bf F}_j(\hat{\mathbf{x}})\\
&\hspace*{3.1cm}-\partial_{\hat{x}_2}\iota_j(\psi_{0,i})(\hat{\mathbf{x}})
\partial_{\hat{x}_1}{\bf F}_j(\hat{\mathbf{x}})\Big)
\end{aligned}
\end{equation}
for all $\mathbf{x}={\bf F} _j(\hat{\mathbf{x}})\in\Gamma _j$, 
$\hat{\mathbf{x}}\in[0,1]^2$, see \cite{DHP16} for example for 
the precise derivation. This amounts to the following expression 
of the hypersingular operator in closed form 
\begin{equation}\label{eq:hypersingularonrefernce}
\begin{aligned}
&\langle\mathcal{W}\psi_{0,k},\psi_{0,\ell}\rangle_{\Gamma}
=
\sum_{i,j=1}^n\int_{[0,1]^2}\int_{[0,1]^2}k_{j,i}(\hat{\bs x},\hat{\bs y})\\
&\quad\times\bigg(
\nabla_{\hat{\bs x}}\iota_j(\psi_{0,k})(\hat{\bs x})^\intercal
K_{j,i}(\hat{\bs x},\hat{\bs y})^{-1}
\nabla_{\hat{\bs y}}\iota_i(\psi_{0,\ell})(\hat{\bs y})
\\
&\quad-\kappa^2
\langle{\bs n}_{\bs x},{\bs n}_{\bs y}\rangle_{\mathbb{R}^3}\iota_j(\psi_{0,k})
(\hat{\bs x})\iota_i(\psi_{0,\ell})(\hat{\bs y})\bigg)\d\hat{\bs y}\d\hat{\bs x},
\end{aligned}
\end{equation}
where the pull-back of the kernel $k_{j,i}$ is given by
\[
k_{j,i}(\hat{\bs x},\hat{\bs y})=k\big({\bf F}_j(\hat{\bs x}),{\bf F}_i(\hat{\bs y})\big)
\]
and $K_{j,i}$ denotes the first fundamental tensor 
of differential geometry,
\[
K_{j,i}(\hat{\bs x},\hat{\bs y})
=
\big[
\langle\partial_{\hat{\bs x}_k}{\bf F}_j(\hat{\bs x}),
\partial_{\hat{\bs x}_l}{\bf F}_i(\hat{\bs y})\rangle_{\mathbb{R}^3}\big]_{k,l=1}^2
\in\mathbb{R}^{2\times 2}.
\]
Compared to the Laplace case, see \cite{DHP16}, we note the 
occurrence of a second term in the regularized representation
\eqref{eq:hypersingularonrefernce}. 
However, this additional term behaves similar to the single 
layer operator and thus poses no further challenges for 
implementation.

For the numerical evaluation of the first term in
\eqref{eq:hypersingularonrefernce}, recall 
that an ansatz function $\psi_{0,j}|_{\Gamma_i}$ on the patch 
$\Gamma _i$ is given by $\psi_{0,j} = \iota_i^{-1}(\hat{\psi})$ for some
$\hat{\psi}\in S_{p,p}(\Xi_{p,m},\Xi_{p,m})$, see \eqref{def::space2}.
There therefore holds
\[
\nabla_{\hat{\bs x}}\iota_i(\psi_{0,i})(\hat{\bs x}) =
\nabla_{\hat{\bs x}}\iota_i\big(\iota_i^{-1}(\hat{\psi})\big)(\hat{\bs x}) =
\nabla_{\hat{\bs x}}\hat{\psi}(\hat{\bs x}).
\]
Thus, for the purposes of implementation, one only has to provide
$S_{p,p}(\Xi_{p,m},\Xi_{p,m})$ and the directional derivatives of 
$S_{p,p}(\Xi_{p,m},\Xi_{p,m})$, which are given by
\begin{align*}
\partial_1S_{p,p}(\Xi_{p,m},\Xi_{p,m})&=S_{p-1,p}(\Xi_{p,m}',\Xi_{p,m}),\\
\partial_2S_{p,p}(\Xi_{p,m},\Xi_{p,m})&=S_{p,p-1}(\Xi_{p,m}',\Xi_{p,m}'),
\end{align*}
where $\Xi_{p,m}'$ denotes the truncation of $\Xi_{p,m}$, i.e., the 
knot vector $\Xi_{p,m}$ without its first and last knot. These 
spline spaces are readily available.

\subsection{Fast multipole method}\label{sec:multipole}
The black-box fast multipole method, cp.\ \cite{Gie01},
relies on a degenerate kernel approximation of the integral
kernel under consideration. Such an approximation is available
in the kernel's far-field, which means that the supports 
of the trial and test functions have to be sufficiently 
distant from each other -- they are \emph{admissible}. 

One arrives at an efficient algorithm, if one 
subdivides the set of trial functions hierarchically 
into so-called clusters. Then, the kernel interaction 
of two clusters is computed by using the degenerate 
kernel approximation if the clusters are admissible.
This means a huge matrix block in the system matrix 
is replaced by a low-rank matrix. If the clusters are not 
admissible, then one subdivides them and considers 
the interactions of the respective children. That way,
the assembly of the Galerkin matrix can be perfomed
in essentially linear complexity. 

For the realization of the multipole method in the present 
context of isogeomtric boundary element methods, we refer 
the reader to \cite{DHK+18,DHP16}. A particular advantage 
of the referred compression method is that the isogeometric 
setting allows to perform the compression of the system 
matrix in the reference domain rather than the computational 
domain. This means that we consider the pull-back of the 
kernel \eqref{eq:pbkernel} instead of the kernel in free space, 
as originally proposed in \cite{HP13}. Thus, the rank of 
the low-rank blocks in the number of one-dimensional 
interpolation points $p$ decreases from $\mathcal{O}(p^3)$ 
to $\mathcal{O}(p^2)$. The compressed matrix is finally 
represented in the \(\mathcal{H}^2\)-matrix format as
usual, see \cite{BH02}.

For the potential evaluation, i.e., for evaluating \eqref{eq:pot1}
and \eqref{eq:pot2}, we exploit a similar approximation
of the kernel function. However, this time we perform the
low-rank approximation in physical space, that is, we 
employ a degenerate kernel approximation for the kernel \(k\). 
Rather than clustering elements as before, we directly cluster
evaluation and quadrature points and realize the potential
evaluations by means
of matrix-vector multiplications. The rank of the low-rank 
blocks is in this case $\mathcal{O}(p^3)$. In particular, we 
may employ a matrix-free version, as all blocks are only 
required once. The advantage of this approach becomes 
immanent if the number of potential evaluation points increases 
proportionally to the number of boundary elements. In this 
case, the cost of the proposed potential evaluations 
scales essentially linearly instead of quadratically. 

\section{Numerical experiments} 	\label{sec:experiments}
\begin{figure}
\includegraphics[width=0.5\textwidth,clip=true,trim=540 290 300 400]{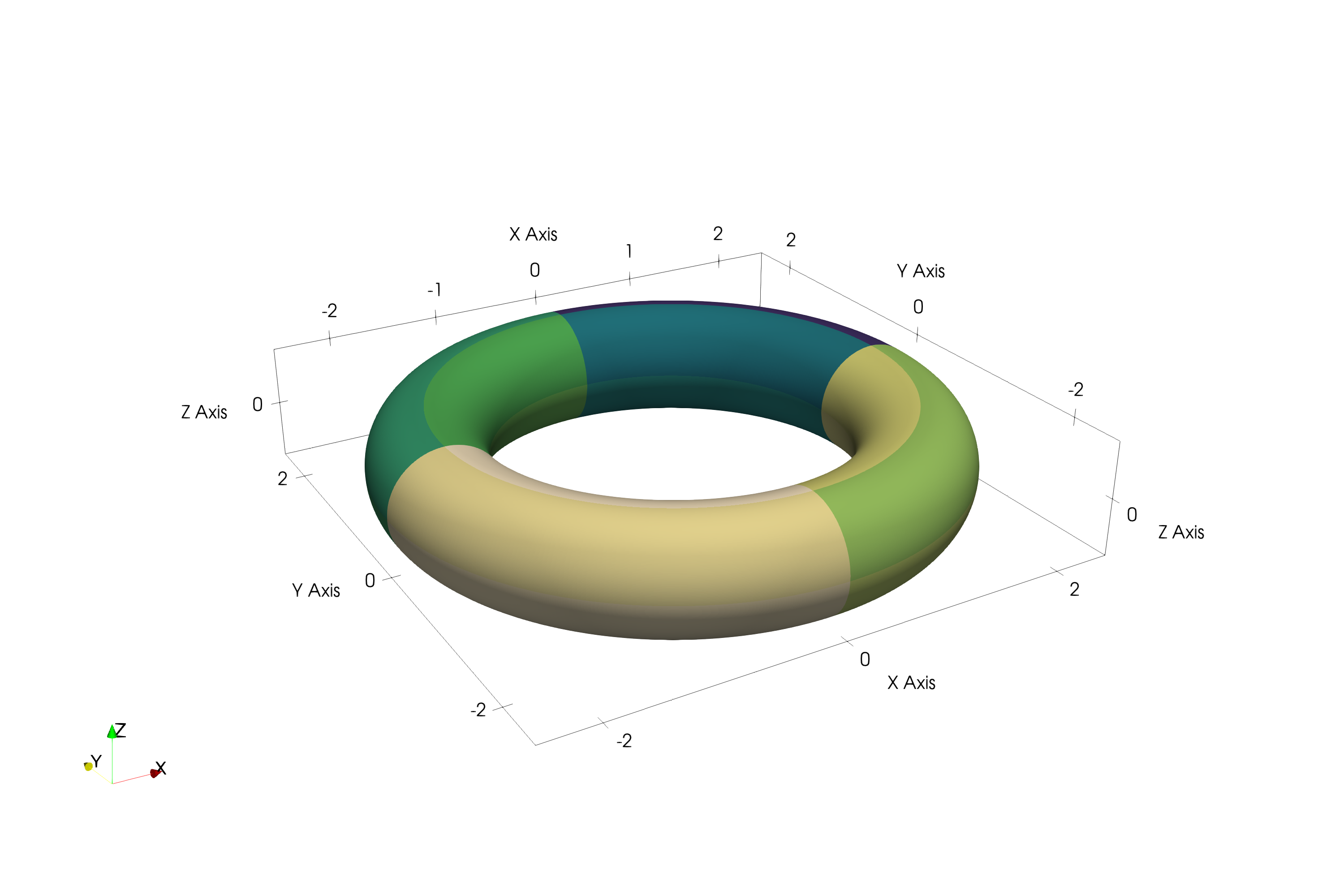}
\caption{\label{fig::torus_setup}Torus represented by 16 patches and
illustration of its dimensions.}
\end{figure}

\subsection{Setup}
The numerical experiments are performed by using the publicly 
available \verb|C++| library \verb|Bembel|, see \cite{bembel,bembelpaper}. To 
this end, the previously not available operators (double layer, adjoint 
double layer, and hypersingular operator) were implemented.
Each of the matrices in the combined field integral equations 
\eqref{eq:CFIED2N} and \eqref{eq:CFIEN2D} was computed
separately in compressed form as $\mathcal{H}^2$-matrix by 
using the fast multipole method on the reference domain from
\cite{bembelpaper,DHP16}. The compression parameters for the fast
multipole method were set to the default values ($\eta=1.6$, nine
interpolation points per direction), see \cite{bembelpaper,DHP16} for
more details. The product of the matrix sums with vectors was 
implemented using lazy evaluation and the arising systems of
linear equations \ref{eq:CFIED2N} and \ref{eq:CFIEN2D} were 
solved up to relative machine precision by means of a restarted 
GMRES method with a restart after 30 iterations. Finally, all 
computations were performed in parallel by using the built-in
\verb|OpenMP|-parallelization of \verb|Bembel| on a compute 
server with 1.3 terrabyte RAM and four Intel(R) Xeon(R) E7-4850 
v2 CPU with twelve 2.30GHz cores each and hyperthreading disabled.

\subsection{Convergence benchmark}
\begin{figure}
\includegraphics[width=0.5\textwidth,clip=true,trim=700 550 400 600]{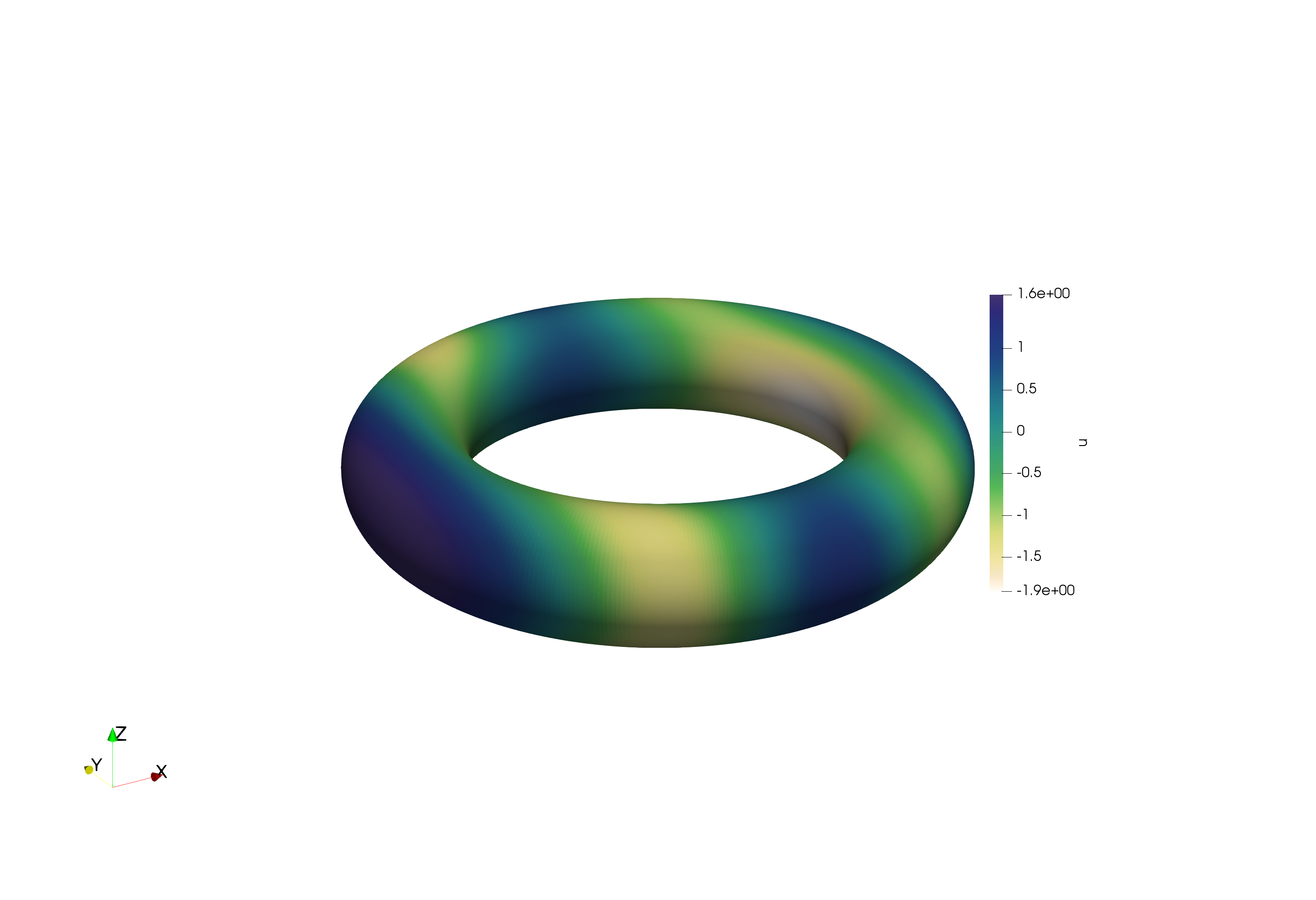}\\
\includegraphics[width=0.5\textwidth,clip=true,trim=700 550 400 600]{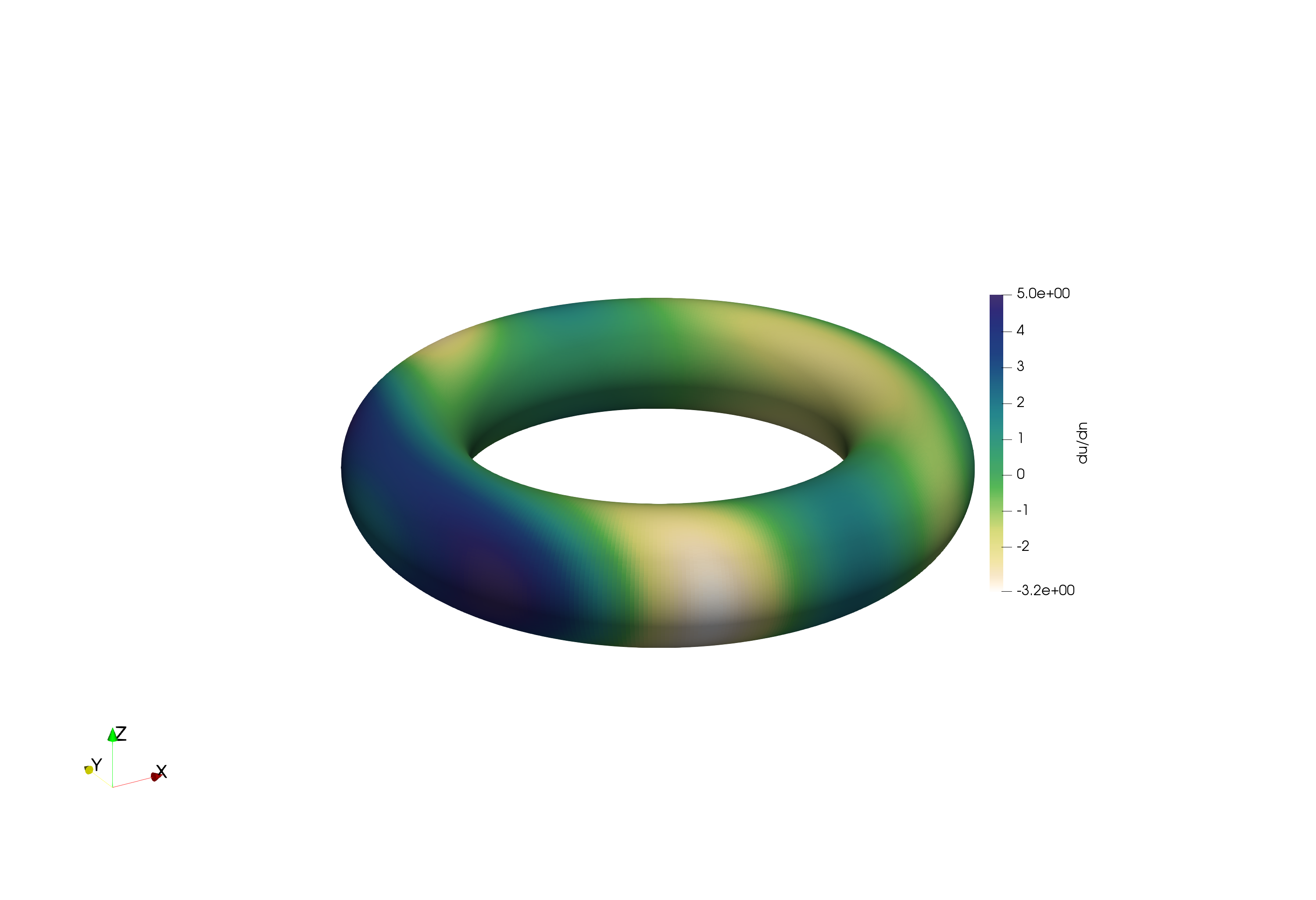}
\caption{\label{fig::torus_wave} The Dirichlet data $u$ (top) of the total wave in 
case of a sound-hard torus and the Neumann data $\partial u/\partial{\bs n}$ 
(bottom) of the total wave in case of a sound-soft torus.}
\end{figure}

In order to study convergence rates, we consider a 
torus with major radius two and minor radius 0.5 that is 
represented by 16 patches, see Figure~\ref{fig::torus_setup}
for an illustration. On this geometry, we aim at computing the 
scattered wave of a plane incident wave in $x$ direction 
with wavenumber 2.5. The scattered wave is then measured
on 100 points distributed on a sphere with radius 5 around the
origin. We refer to Figure~\ref{fig::torus_wave} for an illustration
of the Dirichlet data of the total wave (top plot) in case of a 
sound-hard torus and the Neumann data of the total wave 
(bottom plot) in case of a sound-soft torus.

\begin{figure*}
\begin{tikzpicture}
\begin{semilogyaxis}[
height=0.3\textwidth,
width=0.5\textwidth,
xmin = 0,
xmax = 6,
ymax = 1e1,
ymin = 1e-8,
legend style={legend pos=south west,font=\tiny},
xlabel=refinement level $m$,
ylabel={$\ell^{\infty}$-error},
title={Convergence of sound-soft scattering},
cycle multi list={%
	color list\nextlist
	mark list
},
legend columns=2
]
\foreach \P in {0,...,4} {
	\addplot table[x=M,y=error]{data/log_D2N_\P.log};
	\addlegendentryexpanded{$p=\P$};
}  
\foreach \P in {0,...,4} {
	\addplot[domain=0:6,color=black,dashed] {1.25*pow(4,\P)*pow(2,-2*(\P+1)*x)};
}
\end{semilogyaxis}
\end{tikzpicture}
\begin{tikzpicture}
\begin{semilogyaxis}[
height=0.3\textwidth,
width=0.5\textwidth,
xmin = 0,
xmax = 6,
ymax = 1e1,
ymin = 1e-8,
legend style={legend pos=south west,font=\tiny},
xlabel=refinement level $m$,
ylabel={$\ell^{\infty}$-error},
title={Convergence of sound-hard scattering},
cycle multi list={%
	color list\nextlist
	mark list
},
legend columns=2
]
\addplot coordinates{};
\foreach \P in {1,...,4} {
	\addplot table[x=M,y=error]{data/log_N2D_\P.log};
	\addlegendentryexpanded{$p=\P$};
}  
\foreach \P in {1,...,4} {
	\addplot[domain=0:6,color=black,dashed] {0.3*pow(4,\P)*pow(2,-2*(\P+0.5)*x)};
}
\end{semilogyaxis}
\end{tikzpicture}
\caption{\label{fig:convergence}Convergence of the combined 
field integral equations for various polynomial degrees. The 
dashed lines illustrate the expected convergence rates of 
$\mathcal{O}\big(h^{2p+2}\big)$ in case of sound-soft
obstacles (left) and $\mathcal{O}\big(h^{2p+1}\big)$
in case of sound-hard obstacles (right).}
\end{figure*}
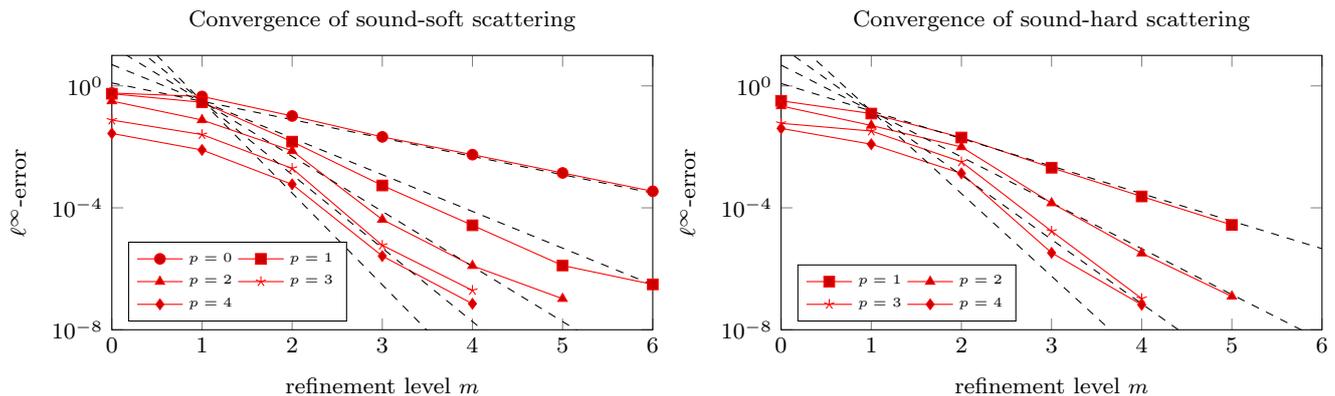

\begin{figure*}
\begin{tikzpicture}
\begin{semilogyaxis}[
height=0.3\textwidth,
width=0.5\textwidth,
xmin = 0,
xmax = 6,
ymin = 1e-1,
ymax = 1e5,
legend style={legend pos=south east,font=\tiny},
xlabel=refinement level $m$,
ylabel={wall clock time (sec.)},
title={Time-to-solution for the sound-soft scattering case},
cycle multi list={%
	color list\nextlist
	mark list
},
legend columns=2
]
\foreach \P in {0,...,4} {
	\addplot table[x=M,y=time]{data/log_D2N_\P.log};
	\addlegendentryexpanded{$p=\P$};
}
\addplot[domain=0:6,color=black,dashed] {x*pow(4,x)};
\end{semilogyaxis}
\end{tikzpicture}
\begin{tikzpicture}
\begin{semilogyaxis}[
height=0.3\textwidth,
width=0.5\textwidth,
xmin = 0,
xmax = 6,
ymin = 1e-1,
ymax = 1e5,
legend style={legend pos=south east,font=\tiny},
xlabel=refinement level $m$,
ylabel={wall clock time (sec.)},
title={Time-to-solution for the sound-hard scattering case},
cycle multi list={%
	color list\nextlist
	mark list
},
legend columns=2
]
\addplot coordinates{};
\foreach \P in {1,...,4} {
	\addplot table[x=M,y=time]{data/log_N2D_\P.log};
	\addlegendentryexpanded{$p=\P$};
}
\addplot[domain=0:6,color=black,dashed] {x*pow(4,x)};

\end{semilogyaxis}
\end{tikzpicture}
\caption{\label{fig:scaling}Scaling of the combined field 
integral equations for various polynomial degrees. The 
dashed lines illustrate log-linear scaling.}
\end{figure*}
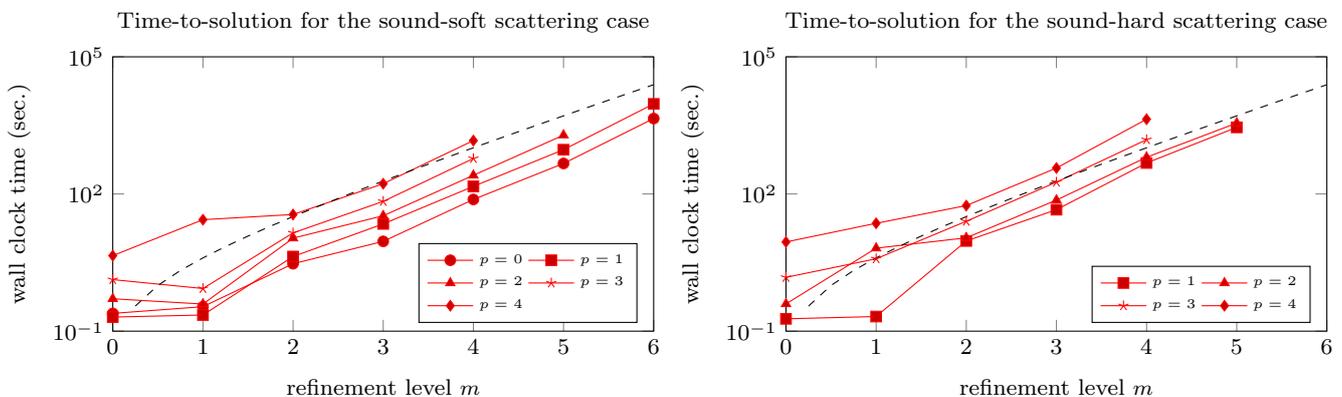

The optimal convergence rates for the potential 
evaluation in case of splines of degree $p$ are 
$\mathcal{O}\big(h^{2p+2}\big)$ for the boundary 
integral equation \eqref{eq:D2N} which corresponds to 
sound-soft obstacles and $\mathcal{O}\big(h^{2p+1}\big)$ 
for the boundary integral equation \eqref{eq:N2D} which 
corresponds to sound-hard obstacles. Since the obstacle
under consideration is smooth, we should achieve these 
convergence rates. Note that these rates are twice as 
high as for the collocation method and are known as 
the \emph{superconvergence} of the Galerkin formulation, 
compare~\cite{Steinbach} for example. Figure 
\ref{fig:convergence} validates that we indeed reach 
these theoretical achievable convergence rates when 
compared to solutions obtained from an indirect 
formulation using a single layer or adjoint double 
layer ansatz, respectively.

Figure \ref{fig:scaling} illustrates the scaling of the runtimes 
of the computations. Instead of a quadratic scaling of the
runtimes, which we would have in the case of a traditional 
boundary element method, one figures out that the 
multipole-accelerated isogeometric boundary element 
method scales essentially linearly as expected. This 
enables large-scale calculations as we will consider 
in the next example.

\subsection{Computational benchmark}
\begin{figure*}
\includegraphics[width=\textwidth,clip=true,trim=500 290 360 400]{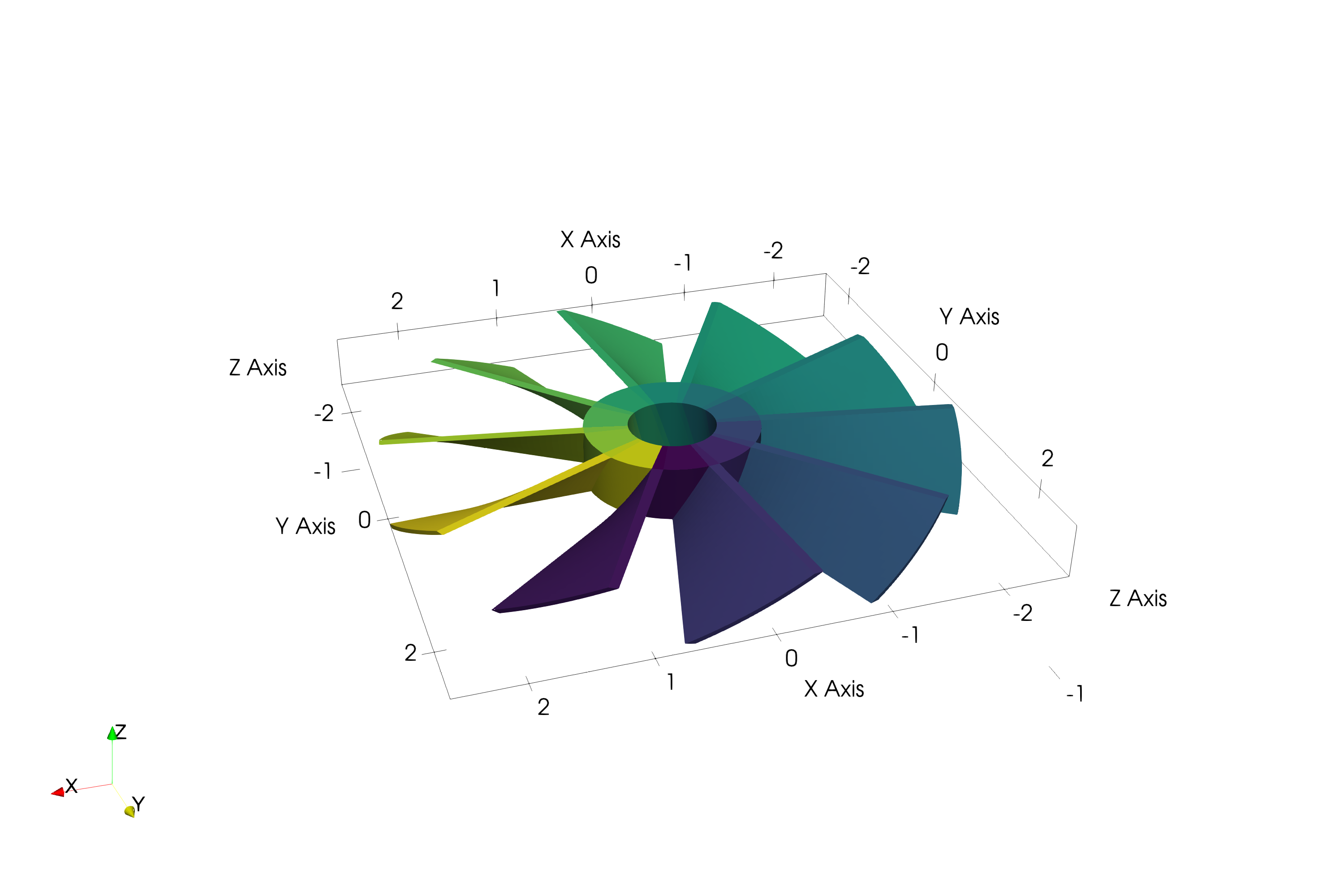}
\caption{\label{fig:turbine}Turbine geometry with 120 patches.}
\end{figure*}

As a computational benchmark, we consider a turbine with 
ten blades that is parametrized by 120 patches as illustrated 
in Figure \ref{fig:turbine}. Thereof, it can be figured out that the
turbine has a diameter of 5. Again, we compute the scattered 
wave of a plane incident wave in $x$ direction, but with 
wavenumber 1.0. 

We choose cubic B-splines and three refinement levels to 
discretize the Cauchy data $u$ and $\partial u/\partial{\bs n}$ 
on the surface geometry. This results in 14'520 degrees of freedom 
in case of a sound-soft turbine and 12'000 degrees of freedom 
in case of a sound-hard turbine, respectively. The overall solution 
time for assembly and solution of the underlying systems of 
linear equations requires only about a few hours.

\begin{figure*}
\includegraphics[height=0.4\textwidth,clip=true,trim=700 600 900 500]{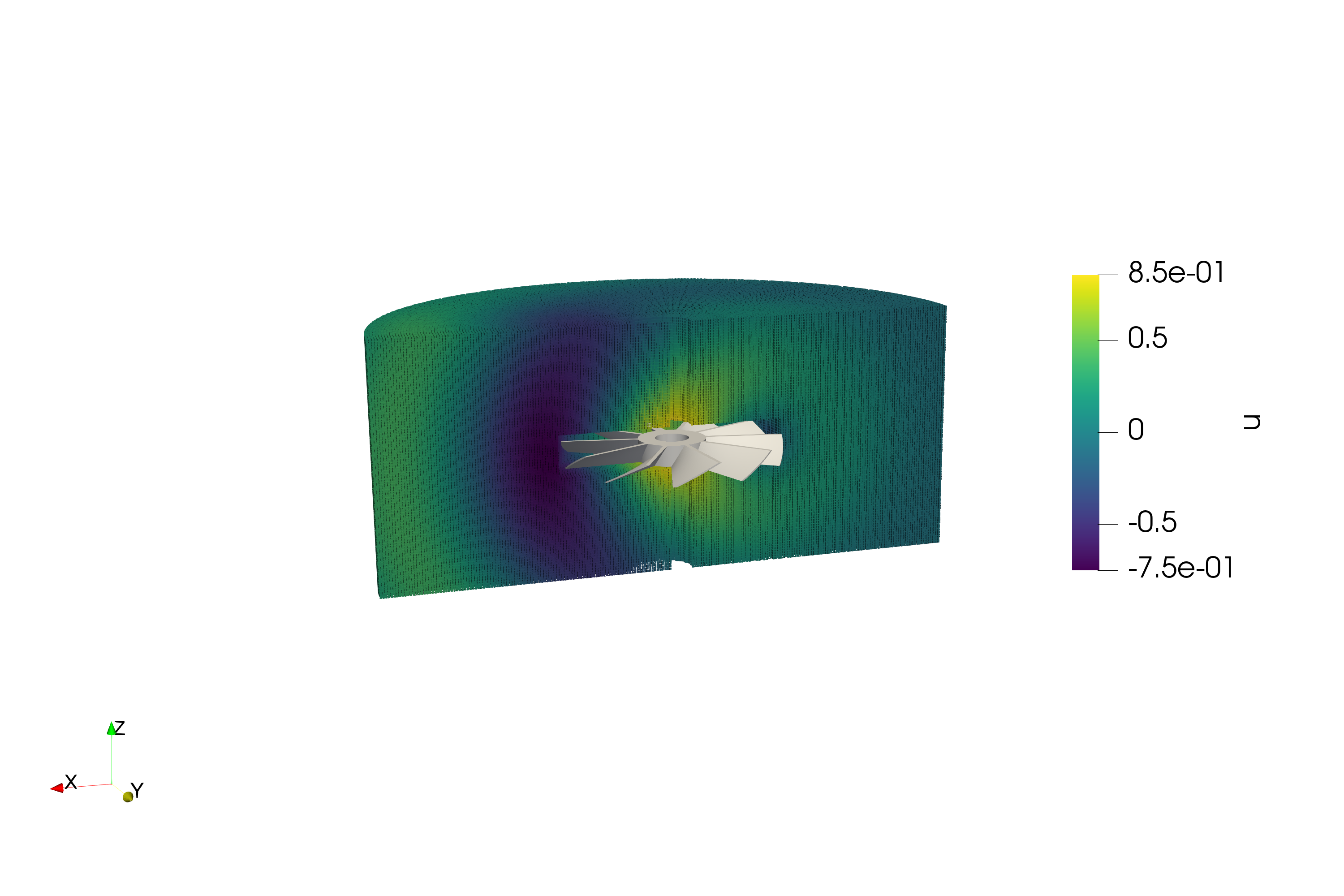}
\includegraphics[height=0.4\textwidth,clip=true,trim=2200 600 100 500]{vtp/turbine/scatter_u}
\caption{\label{fig:scatteredwave}Scattered wave evaluated in 3'664'832 points for the sound-hard case.}
\end{figure*}

\begin{figure*}
\begin{tikzpicture}
\begin{loglogaxis}[
height=0.3\textwidth,
width=0.5\textwidth,
legend style={legend pos=north west,font=\tiny},
xlabel=number of potential points,
ylabel={wall clock time (sec.)},
title={Potential evaluation time, sound-soft scattering},
cycle multi list={%
	color list\nextlist
	mark list
}
]
\addplot table[x=gridpoints,y=tfull]{data/log_D2N_pot_p3_l3.log};
\addlegendentryexpanded{conventional potential evaluation};
\addplot table[x=gridpoints,y=tH2]{data/log_D2N_pot_p3_l3.log};
\addlegendentryexpanded{$\mathcal{H}^2$ potential evaluation};
\addplot[domain=1e3:1e7,color=black,dashed] {1e-3*x};
\addlegendentryexpanded{linear scaling};
\end{loglogaxis}
\end{tikzpicture}
\begin{tikzpicture}
\begin{loglogaxis}[
height=0.3\textwidth,
width=0.5\textwidth,
legend style={legend pos=north west,font=\tiny},
xlabel=number of potential points,
ylabel={wall clock time (sec.)},
title={Potential evaluation time, sound-hard scattering},
cycle multi list={%
	color list\nextlist
	mark list
}
]
\addplot table[x=gridpoints,y=tfull]{data/log_N2D_pot_p3_l3.log};
\addlegendentryexpanded{conventional potential evaluation};
\addplot table[x=gridpoints,y=tH2]{data/log_N2D_pot_p3_l3.log};
\addlegendentryexpanded{$\mathcal{H}^2$ potential evaluation};
\addplot[domain=1e3:1e7,color=black,dashed] {1e-3*x};
\addlegendentryexpanded{linear scaling};
\end{loglogaxis}
\end{tikzpicture}
\caption{\label{fig::potentialscaling}Computation time of
conventional and $\mathcal{H}^2$-matrix
accelerated potential evaluation for various number of points.}
\end{figure*}
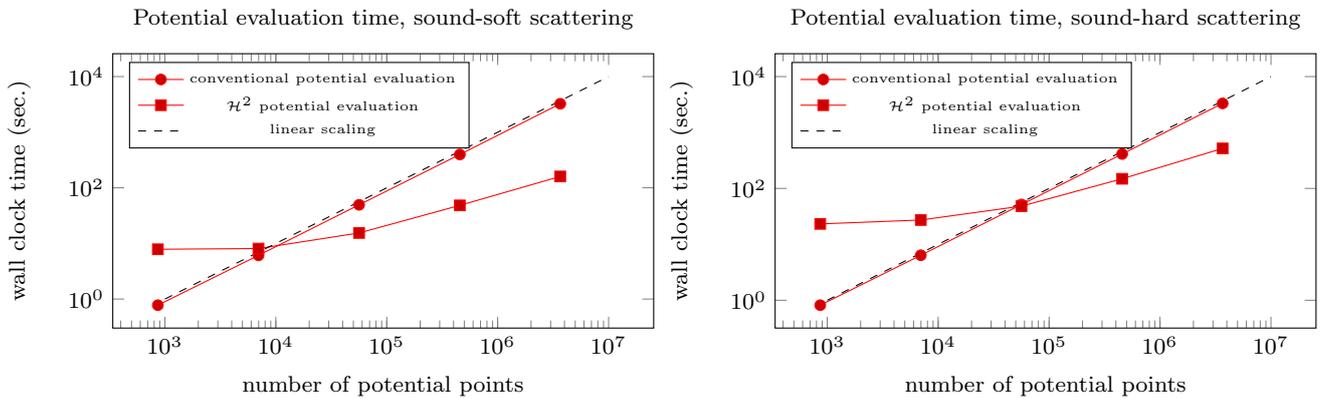

We compute next the scattered wave in a cylinder on up 
to 3'664'832 points, compare Figure~\ref{fig:scatteredwave} 
for an illustration. To demonstrate the efficiency of the 
fast potential evaluation, we compare the scaling of the 
multipole-accelerated potential evaluations with the 
traditional potential evaluations. Figure 
\ref{fig::potentialscaling} illustrates that -- after a 
certain warm-up phase for only a few potential points -- 
the $\mathcal{H}^2$-matrix accelerated potential 
evaluation is indeed superior to the conventional one
when increasing the number of evaluation points.
Consequently, the calculation of the scattered wave
also in free space becomes feasible and very efficient.

\section{Conclusion} 	\label{sec:conclusion}
We have solved acoustic scattering problems for sound-soft and
sound-hard scatterers by means of frequency stable combined field
integral equations and the isogeometric boundary
integral equation method. The major
advantage of this approach is that no mesh for the unbounded 
exterior domain is required. Our discretization
is based on Galerkin's method together with an appropriate
regularization of the hypersingular operator. The method becomes
computationally efficient by the use of a black-box fast
multipole method tailored to isogeometric surfaces. A similar
approach is employed for the postprocessing step to efficiently 
evaluate the solution in free space. We have presented convergence 
benchmarks that impressively demonstrate the high accuracy of 
the isogeometric boundary element method. In addition, we have 
considered a complex computational benchmark on a complex 
geometry, which corroborates the feasibility of the
approach in the engineering practice.
 
\bibliographystyle{plain}

\end{document}